\documentclass[a4paper, mathscr, 12 pt]{article}
\usepackage{eucal, amssymb,amsmath,graphicx, amsfonts,amsthm, latexsym, enumerate}
\def\bold1{\boldsymbol{1}}
\def\bold0{\boldsymbol{0}}

\numberwithin{equation}{section}

\begin{document}

\title{Inhomogeneous epidemics on weighted networks}

\author{ Tom Britton, Stockholm University\thanks{Department of
Mathematics, Stockholm University, SE-106 91 Stockholm, Sweden. {\it
E-mail}: tom.britton@math.su.se} \\
David Lindenstrand, Stockholm University\thanks{Department of
Mathematics, Stockholm University, SE-106 91 Stockholm, Sweden. {\it
E-mail}: davlin@math.su.se}\ \thanks{\emph{To whom correspondence
should be addressed.}}}
\date{\today}
\maketitle

\begin{abstract}

\noindent A social (sexual) network is modeled by an extension of the configuration model to the situation where edges have weights, e.g.\ reflecting the number of sex-contacts between the individuals. An epidemic model is defined on the network such that individuals are heterogeneous in terms of how susceptible and infectious they are. The basic reproduction number $R_0$ is derived and studied for various examples, but also the size and probability of a major outbreak. The qualitative conclusion is that $R_0$ gets larger as the community becomes more heterogeneous but that different heterogeneities (degree distribution, weight, susceptibility and infectivity) can sometimes have the cumulative effect of \emph{homogenizing} the community, thus making $R_0$ smaller. The effect on the probability and final size of an outbreak is more complicated.
\vskip1cm
\end{abstract}

Keywords: basic reproduction number, heterogeneity, random network, stochastic epidemic model.


\section{Introduction}

\noindent Epidemic models have a long history in mathematical modelling (see e.g.\ Diekmann and Heesterbeek \cite{DH00} for an overview). Early models assumed a homogeneous community but later this was relaxed be allowing individuals to vary, for example by dividing individuals into different groups, thus defining so-called multitype epidemics. More recent models admitting local structures in the community have been included into epidemic models, household models and network models being the two main examples (e.g.\ Ball et al.\ \cite{BMS97} and Andersson \cite{A99}). Admitting local structure have the effect that stochastic models are favourable in that when only few individuals affect the risk of becoming infected the outcome should be random.\\

\noindent The current paper aims at combining the two types of heterogeneities mentioned above: individual heterogeneities with network models. More precisely, we want to analyse how individual variation in susceptibility and infectivity affect the epidemic spread in a population composed into a social network with weighted edges. We have a sexually transmitted infection (STI) as a motivation for the paper. This heterogeneity between individuals might then correspond to varying sexual risk behaviour, e.g.\ not using a condom, or physiological features such as more or less susceptible Mucous membranes. The weights on the edges can for example correspond to the number of sexual contacts between the two individuals in question. To this end we extend the work in \cite{BD11}, which studies a weighted configuration model, which describes a stochastic epidemic model on a network, by including heterogeneity in susceptibility and infectivity, i.e. the ability to receive, and transmit, the infection.  For a fixed individual $i$, we let the susceptibility and infectivity be an outcome of a pair of random variables $(X_i, Y_i)$ where these two random variables may depend on each other, but are assumed independent between individuals.\\

\noindent We derive the basic reproduction $R_0$ and investigate how $R_0$ is influenced by the coefficients of variation $CV_X=\sigma_X/\mu_X$ and $CV_Y=\sigma_Y/\mu_Y$, for different correlations, $\rho_{X,Y}$, between $X$ and $Y$. In particular we compare the result to the $R_0$ obtained for fixed susceptibility and infectivity for all individuals ($CV_X=CV_Y=0$), corresponding the model analysed in \cite{BD11}. We also investigate how the probability and size of a large outbreak depends on $CV_X$ and $CV_Y$.\\

\noindent For an introduction to network models and their applications (including epidemics) we refer to Newman \cite{newman}, whereas van der Hofstad \cite{hofstad} gives a more technical and exhaustive treaties of the random graphs and their properties, e.g.\ thorough treatment of the configuration model, and the relation between branching processes and random graphs. Epidemics on networks allowing individual heterogeneities have been studied earlier, but has perhaps not yet received enough attention. Miller \cite{miller} studies an epidemic on an unweighted graph where susceptibility and infectivity varies among the individuals, and derives bounds on the probability and size of a large outbreak. A similar problem is studied by Trapman and Meester in \cite{trapman} who define a percolation model and derive bounds for e.g.\ expected final size and outbreak probability as function of the infectivity and susceptibility.\\

\noindent The rest of the paper is outlined as follows. In Section \ref{secmod} we specify the graph model, the epidemic model and give a general expression for the basic reproduction number $R_0$. In Section \ref{secr0} we investigate the effect of individual heterogeneity on $R_0$ for different examples of the model. In Section \ref{secpi} we derive the outbreak probability $\pi$ in a simple case and in Section \ref{secdis} we conclude with a short discussion.

\section{Model}\label{secmod}

\subsection{A model for weighted network}

The network model we study was originally defined by Britton et al.~\cite{BD11}, which is an extension of the configuration model (see e.g. \cite{hofstad}) but allowing for edges to have different weights. It is defined as follows. We have a population of size $n$. Individual $i$ has a random number $D_i$ of half-edges where $\{D_i\}$ are iid\ with probability function $\{p_D(k), k \geq 0 \}$ ($D_i$ is called the degree of individual $i$).  Furthermore, each half-edge of individual $i$ has a random weight $W_{ij}\in\mathbb{N}$, $j=1,\dots ,D_i$, possibly dependent of $D_i$ but being independent of each other (in the STI application $W_{ij}$ can for example reflect the number of sexual contacts individual $i$ has with its partner $j$). The probability that a random half-edge of an individual with $d$ half-edges has weight $w$, is denoted $q(w|d)$. The network is formed by randomly connecting half-edges with the same weight. If the number of half-edges of a specific weight is odd, the last half-edge is ignored. We assume that the variance of $D$ is finite. Thus the number of self loops and multiple edges are negligible if $n$ is large (see Section 7 in \cite{hofstad}).\\

\subsection{An epidemic model on the weighted network}

\noindent To each individual $i$ we assign a random vector $(X_i,Y_i)$ with probability function $p_{X,Y}$ being independent of $(D_i, W_{i1}, \dots ,W_{iD_i})$. The first component $X_i$ denotes the susceptibility and $Y_i$ the infectivity (in case of getting infected) of individual $i$. We assume that $0\leq X_i\leq 1$ and $0\leq Y_i \leq 1$ and allow the variables $X_i$ and $Y_i$ to be dependent. An individual $i$ is hence described by its individual susceptibility and infectivity $(X_i,Y_i)$, together with the independent random vector $(D_i, W_{i1}, \dots ,W_{iD_i})$ indicating how many neighbours $i$ has and the weights on the edges connecting to them.\\

\noindent Consider two individuals directly connected to each other by an edge of weight $w$. Suppose the first individual has susceptibility/infectivity $(x_1,y_1)$ and that he/she is infected by someone else, and suppose the second individual is still susceptible having susceptibility/infectivity $(x_2,y_2)$. The epidemic model is then defined by saying that the first individual infects the second with probability $t(w,y_1,x_2)$ defined by
\begin{equation}\label{eqprob}
 t(w,y_1, x_2)=1-(1-y_1x_2)^w.
\end{equation}
The intuition behind Equation (\ref{eqprob}) is the following. In one contact, the first individual infects the second with probability $x_2y_1$ (the more infectious the first is and the more susceptible the second is the higher risk of disease transmission). Hence the first individual does \textit{not} infect the second individual in $w$ contacts with probability $(1-x_2y_1)^w$. Given all susceptibilities, infectivities, degrees and weights, transmission events are defined to be mutually independent. Initially one randomly selected individual is externally infected (the index case) and the rest of the community is susceptible to the disease. The index case infects a random number of his/her neighbours following the transmission probability defined in (\ref{eqprob}) and then becomes immune. These newly infected may in turn infect some of their not yet infected members and then become immune. The epidemic continues until there are no new infections. Then the epidemic stops. Those who were infected make up the \emph{final outcome} of the epidemic and the number of infected is called the \emph{final size} of the epidemic.\\

\noindent We now assume that the size $n$ of the community is large. The initial stage (before a non-negligible fraction have been infected) can then be approximated by a multiptype branching process (e.g.\ \cite{J75}) as is nearly always the case with epidemics in large populations (cf.~\cite{BD11} for the case with unweighted edges). In the current model the \emph{type} of the individual in the branching process approximation is characterised by the degree of the individual together with the susceptibility and infectivity. In order to compute the basic reproduction number $R_0$, loosely defined as the expected number of new infections caused by a random infected during the early stages of the epidemic, we first compute the expected number of individuals of a given type that an infected of a given type infects during the early stages of the epidemic. In branching process terminology this is the mean offspring matrix.\\

\noindent Let $m_{(d_1,x_1,y_1),(d_2,x_2,y_2)}$ be the mean number of individuals of type $(d_2,x_2,y_2)$ that gets infected by one $(d_1,x_1,y_1)$-individual during the early stages of the epidemic, and let $M=\{m_{(d_1,x_1,y_1),(d_2,x_2,y_2)}\}$ be the mean offspring matrix. We now derive an expression for $m_{(d_1,x_1,y_1),(d_2,x_2,y_2)}$ in a similar way as in \cite{BD11}. Let $(d_1,x_1,y_1)$ be a fixed infected individual in the early stage of the epidemic outbreak, and let $p_{(d_1,x_1,y_1)}(d_2,x_2,y_2)$ denote the probability that it infects a $(d_2,x_2,y_2)$-individual along one of its $d_1-1$ edges to susceptibles (the individual was infected through one of its edges). Furthermore, let $Z_{(d_1,x_1,y_1)}(d_2,x_2,y_2)$ denote the (random) number of individuals of type $(d_2,x_2,y_2)$ it infects; since type and infections along different edges are independent it follows that $Z_{(d_1,x_1,y_1)}(d_2,x_2,y_2)$ is binomially distributed with parameters $d_1-1$ (in the early stages all but its ''infector'' are susceptible) and $p_{(d_1,x_1,y_1)}(d_2,x_2,y_2)$. From this it follows that $m_{(d_1,x_1,y_1),(d_2,x_2,y_2)}=(d-1)p_{(d_1,x_1,y_1)}(d_2,x_2,y_2)$. It remains to compute $p_{(d_1,x_1,y_1)}(d_2,x_2,y_2)$. \\

\noindent We compute $p_{(d_1,x_1,y_1)}(d_2,x_2,y_2)$ by summing over all possible weights along the edge in question, since the transmission probability but also the type of the connected node depends on the weight. For a given weight $w$ the probability $\tilde p_w(d_2, x_2, y_2)$ that the edge connects to a $(d_2,x_2,y_2)$-individual is proportional to $d_2p_D(d_2)$ (since the fraction of edges connecting to $d_2$-individuals is $d_2p_D(d_2)/\sum_d dp_D(d)$), and proportional to $q(w|d_2)$, the latter being the probability that an individual of degree $d_2$ has weight $w$ along a given edge. The probability $\tilde p_w(d_2, x_2, y_2)$ is also proportional to $p_{X,Y}(x_2,y_2)$ since the susceptibility and infectivity are independent of the network structure. As a consequence we have
\begin{equation}\label{m1}
 \tilde p_{w}(d_2,x_2,y_2) = \frac{q(w|d_2)d_2p_D(d_2)p_{X,Y}(x_2,y_2)}{\sum_d q(w|d)dp_D(d)}.
\end{equation}

\noindent The quantity $\tilde p_{w}(d_2,x_2,y_2)$ is the probability that an edge with weight $w$ connects to a $(d_2,x_2,y_2)$-individual. We want to compute $p_{(d_1,x_1,y_1)}(d_2,x_2,y_2)$, the probability that a $(d_1,x_1,y_1)$-individual infects a  $(d_2,x_2,y_2)$-individual along a given edge. This probability is obtained by summing over all possible weights, and using  $\tilde p_{w}(d_2,x_2,y_2)$ defined in (\ref{m1}), but multiplied by the probability that the first individual has this weight $w$ and multiplied by the transmission probability $t (w,x_2, y_1)$. We hence get
\begin{equation}\label{m2}
  p_{(d_1,x_1,y_1)}(d_2,x_2,y_2)=\sum_{w}q(w|d_1)t(w,x_2,y_1)\tilde p_{w}(d_2,x_2,y_2).
\end{equation}
We are now ready to compute the mean offspring matrix $M$ with elements $m_{(d_1,x_1,y_1),(d_2,x_2,y_2)}$ denoting the \emph{expected} number of $(d_2,x_2,y_2)$-individuals that one infected $(d_1,x_1,y_1)$-individual infects during the early stages of an outbreak. The corresponding \emph{random} number is $Z_{(d_1,x_1,y_1)}(d_2,x_2,y_2)$, and from before we know that this quantity is binomially distributed with parameters $d_1-1$ and $p_{(d_1,x_1,y_1)}(d_2,x_2,y_2)$ (defined in \ref{m2}). It hence follows that
\begin{eqnarray}\label{m0}
m_{(d_1,x_1,y_1),(d_2,x_2,y_2)}  &=& E[Z_{(d_1,x_1,y_1)}(d_2,x_2,y_2)] \nonumber \\
                             &=& (d_1-1) p_{(d_1,x_1,y_1)}(d_2,x_2,y_2)\nonumber \\
                             &=& (d_1-1)\sum_{w}t(w,x_2,y_1)q(w|d)\tilde p_{w}(d_2,x_2,y_2).
\end{eqnarray}
The basic reproduction number, $R_0$,  is the largest eigenvalue of the matrix $M$. This quantity plays a key role in epidemics and branching processes in that it is a threshold parameter. More precisely, a large epidemic outbreak is possible if and only if $R_0>1$ (see e.g.\ \cite{J75}). Unfortunately it is not possible to state any general features of $R_0$ for the general case, except that it is the largest eigenvalue of $M$. In the next section we study some specific examples.


\section{$R_0$ for various weight, degree, susceptibility and infectivity distributions}\label{secr0}

In order to gain insight in how $R_0$ depends on different heterogeneities: the degree distribution, the weights, and in particular the variable susceptibility and infectivity, we now study a few examples analytically and/or numerically.

\subsection{Unweighted network}\label{unweight}

We begin by analysing the case with fixed weight $W\equiv1$, i.e.\ an unweighted network (this model is a special case of the model analysed in \cite{trapman}). Thus we get the original configuration model as our social network, and the epidemic model where there is heterogeneity in terms of susceptibility and infectivity. In this case Equation (\ref{m1}) simplifies to
\begin{eqnarray*}
 \tilde p_{w}(d_2,x_2,y_2) 
                     &=&  \frac{d_2p(d_2)p_{X,Y}(x_2,y_2)}{\mu_D},
\end{eqnarray*}
where $\mu_D=\sum_ddp_D(d)$ is the mean degree.\\

\noindent It follows that Equation (\ref{m0}) simplifies to \begin{eqnarray}\label{eqR_0}
m_{(d_1,x_1,y_1),(d_2,x_2,y_2)} &=& (d_1-1) p_{(d_1,x_1,y_1)}(d_2,x_2,y_2)\nonumber\\
                                &=&
                                (d_1-1) y_1 \frac{x_2d_2p_D(d_2)p_{X,Y}(x_2,y_2)}{\mu_D}.
\end{eqnarray}
From (\ref{eqR_0}) we see that the elements of $M$ can be written as a product of two factors, one depending on $(d_1,x_1,y_1)$ and the other depending on $(d_2,x_2,y_2)$.
The basic reproduction number $R_0$, i.e. the largest eigenvalue of $M$, is derived as follows. Let $\{ \lambda_i\}$ be the eigenvalues of $M$. Then $trace(M)=\sum_i \lambda_i$. We see that there exists vectors $a$ and $b$ (functions of $(d_1,x_1,y_1)$, $(d_2,x_2,y_2)$) such that $M=a*b^T$ and $Ma=trace(M)a$. Thus $a$ is a eigenvector with eigenvalue $trace(M)$. It follows that all other eigenvalues must be zero and $trace(M)$ hence equals the largest one (see e.g.\ \cite{BM90}). To conclude, we have \\
\begin{eqnarray}
 R_0 &=& \sum_{d,x,y} m_{(d,x,y),(d,x,y)}\label{eqr01}\\
     &=& \sum_{d,x,y} (d-1) x y \frac{d p_D(d)p_{X,Y}(x,y)}{\mu_D}\nonumber\\
     &=& E(XY)(E(D(D-1))\nonumber\\
     &=& \mu_X\mu_Y(1+CV_X CV_Y\rho_{X,Y})\Big(\mu_D+\frac{\sigma^2_D-\mu_D}{\mu_D}\Big).\label{R_0-fixed}
\end{eqnarray}

\noindent We have assumed that $X$ and $Y$ are both discrete. If they instead were continuous then $R_0$ would be the largest eigenvalue of a related functional, but since we can approximate a continuous distribution arbitrarily well by a discrete distribution, we would get the same expression, so (\ref{R_0-fixed}) applies whenever the network is unweighted. \\

\noindent From (\ref{R_0-fixed}) we see that, for an unweighted network, $R_0$ increases with the coefficient of variation of the susceptibility and infectivity for the more likely scenario that infectivity and susceptibility are positively correlated. If on the other hand $\rho_{X,Y}<0$ then $R_0$ is \emph{decreasing} in the coefficients of variation.\\

\begin{figure}[!h]
\centering
\includegraphics[width=7cm, height=7cm]{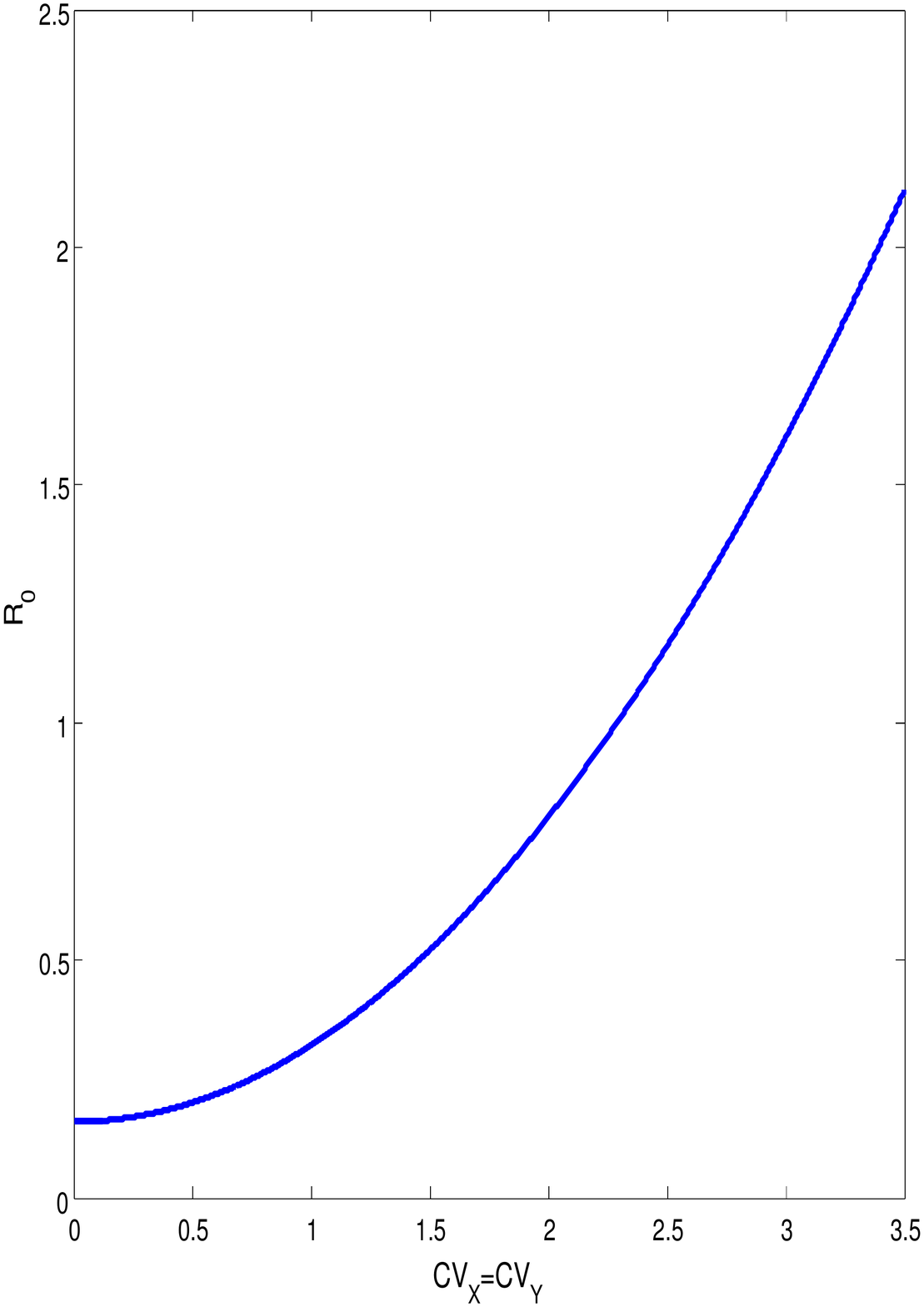}
\caption{$R_0$ as function of $CV_X=CV_Y$ in the case with fixed weight $W\equiv 1$ and degree $D\equiv 5$, and where $\mu_X=\mu_y=0.2$ and $X$ and $Y$ have correlation $0.7$.}
\label{fig12}
\end{figure}

\noindent In Figure \ref{fig12} this situation is illustrated for the case $D\equiv 5$, $\rho_{X,Y}=0.7$ and $\mu_x=\mu_y=0.2$, and assuming the same coefficient of variation in infectivity and susceptibility, i.e.\ that $CV_X=CV_Y$. We see that $R_0$ is increasing with $CV_X=CV_Y$, which is a measure of the heterogeneity in the population. Larger $CV_X=CV_Y$, together with positive correlation, results in more individuals with high ability to transmit $\textit{and}$ receive disease, and more individuals with low ability to transmit \textit{and} receive disease. In Figure \ref{fig12} it is seen that $R_0$ increases with $CV_X=CV_Y$ in this case. Note that the distribution of $(X,Y)$ need not be fully specified, only the mean and coefficient of variation of $X$ and $Y$, together with the mutual correlation is needed.\\

\noindent For the special case that $X \equiv Y$, i.e. the infectivity and susceptibility are identical (fully correlated, $\rho_{X,Y}=1$) but different between individuals, we get
\begin{equation}
 R_0 = \mu_X^2(1+CV_X^2)\Big(\mu_D+\frac{\sigma^2_D-\mu_D}{\mu_D}\Big)\label{r0unw2}.
\end{equation}
From (\ref{r0unw2}) we see that, in the situation where susceptibility and infectivity are identical, $R_0$ increases in the randomness of the susceptibility/infectivity.

\subsection{$d$-networks with random weights but fixed infectivity and susceptibility}\label{secwe}

Suppose now that $D\equiv d$ and $X$ and $Y$ are deterministic (these are both set to $ \sqrt{p}$, so that the probability of infection in one contact becomes $p$). The weight $W$ is random with distribution $q(w|d)=q(w)$, iid\ between all different edges. This is a special case of the model introduced in \cite{BD11}. Since there is only one type of individual in this case, we get that
\begin{eqnarray}
R_0 &=& (d-1)\sum_{w}t(w, \sqrt{p}, \sqrt{p})q(w)\nonumber \\
    &=& (d-1)\sum_{w}(1-(1-p)^w)q(w)\nonumber \\
    &=& (d-1)(1 - G(1-p)),\label{R_0-wei}
\end{eqnarray}
where $G$ is the probability generating function of the weight $W$. Since $t(w, \sqrt{p}, \sqrt{p})=1-(1-p)^w$ is a concave function, we get the from Jensen's inequality that $R_0=E(\pi(W, \sqrt{p}, \sqrt{p}))<\pi(E(W), \sqrt{p}, \sqrt{p})$. Thus, $R_0$ is larger when all edges have the same weight compared to the case with random weights.\\

\begin{figure}[!h]
\centering
\includegraphics[width=7cm, height=7cm]{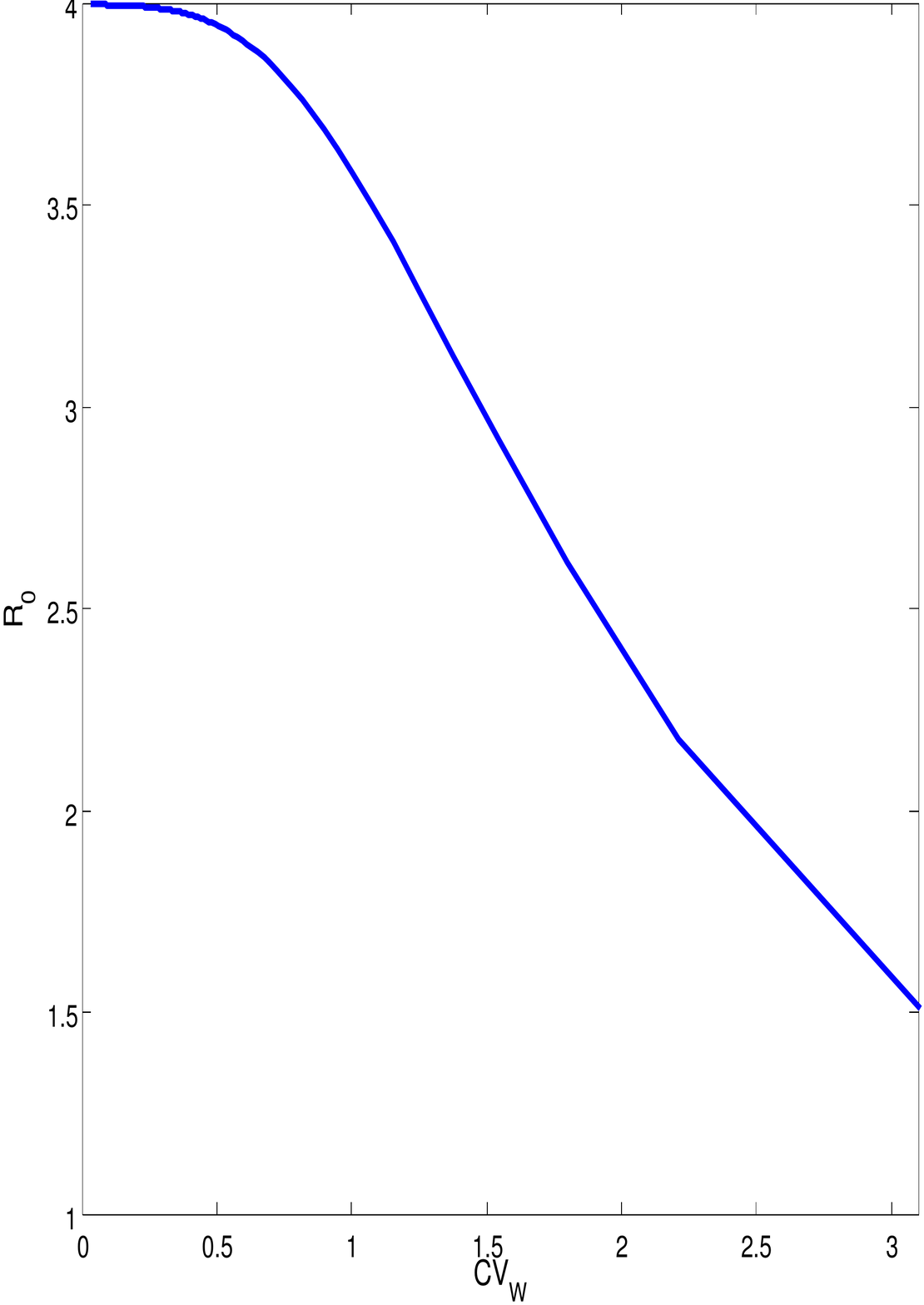}
 \caption{$R_0$ as function of $CV_W$ where $D\equiv 5$, $\mu_W=10$ and $X \equiv Y \equiv \sqrt{0.5}$, i.e. the probability of transmission in one contact is $p=0.5$.}
\label{fig6}
\end{figure}

\noindent In order to obtain an explicit expression we consider the specific example where $W$ follows a negative binomial distribution $W\sim NB(r,\phi)$, i.e.\ where $P(W=w)=\binom{w-1}{r-1}\phi^r(1-\phi)^{w-r}$, for $w=r, r+1,\dots$. The mean, variance and coefficient of variation are given by
\begin{equation}
E(W)=\mu_W=\frac{r}{\phi},\qquad V(W)=\frac{r(1-\phi)}{\phi^2},\qquad CV_W=\sqrt{\frac{1}{r}-\frac{1}{\mu_W}}.\label{mom-negbin}
\end{equation}
The probability generating function equals $G(s)=\left( s\phi/(1-s(1-\phi))\right)^r$. We want to study how $R_0$ is affected by the randomness in the weight distribution. To this end we fix the mean $\mu_W$ and hence set $\phi=r/\mu_W$. Having fixed $\mu_W$ we see from (\ref{mom-negbin}) that the coefficient of variation is decreasing in the remaining parameter $r$. Inserting the probability generating function for the negative binomial distribution into (\ref{R_0-wei}) gives us
\begin{equation}
R_0   = (d-1) \Big(1-\Big(\frac{(1-p)\frac{r}{\mu_W}}{(1-p)\frac{r}{\mu_W}+p}\Big)^{r}\Big).\label{r034}
\end{equation}
 Having the transmission probability $p$ and the expected weight $\mu_W$ fixed we study how $R_0$ depends on the randomness in $W$ by how $R_0$ depends on $r$. Since $R_0$ is an increasing function of $r$ (see section \ref{secapp}), and $CV_W$ is a decreasing function of $r$, we conclude that $R_0$ is a decreasing function of $CV_W$ (and in $V(W)$). This situation is illustrated in Figure \ref{fig6} with $\mu_W=10$, $D\equiv 5$ and $p=0.5$, where this decay is confirmed.

\subsection{Further examples}\label{secneg}

\noindent In the subsections above we have given two special cases of the general model where $R_0$ can be computed explicitly. In general, i.e.\ where weights, degrees, susceptibility and infectivity are all random,  this is not the case. Below we illustrate a few ''toy examples'' having all four features random, with the aim of illustrating different qualitative aspects of how $R_0$ depends on the randomness of the susceptibility and infectivity assuming these are positively correlated. In the examples we assume that the community proportion with high susceptibility is equal to the community proportion with low susceptibility and the same is assumed regarding infectivity. Varying the coefficient of variation in susceptibility/infectivity is achieved by varying the difference between high and low susceptibility/infectivity. Furthermore, the correlation between susceptibility and infectivity is managed by altering the proportion with low susceptibility \textit{and} high infectivity (or equivalently the proportion with high susceptibility \textit{and} low infectivity). In Example 1-3 the distribution of susceptibility and infectivity $(X,Y)$, and the degree distribution, is the same, but with different weight distribution. In Example 1 and 3, $R_0$ is illustrated as a function of $CV_X=CV_Y$, whereas in Example 2 is illustrated as a function of $CV_X$ with $CV_Y$ fixed.\\

\begin{figure}[!h]
\centering
\includegraphics[width=7cm, height=7cm]{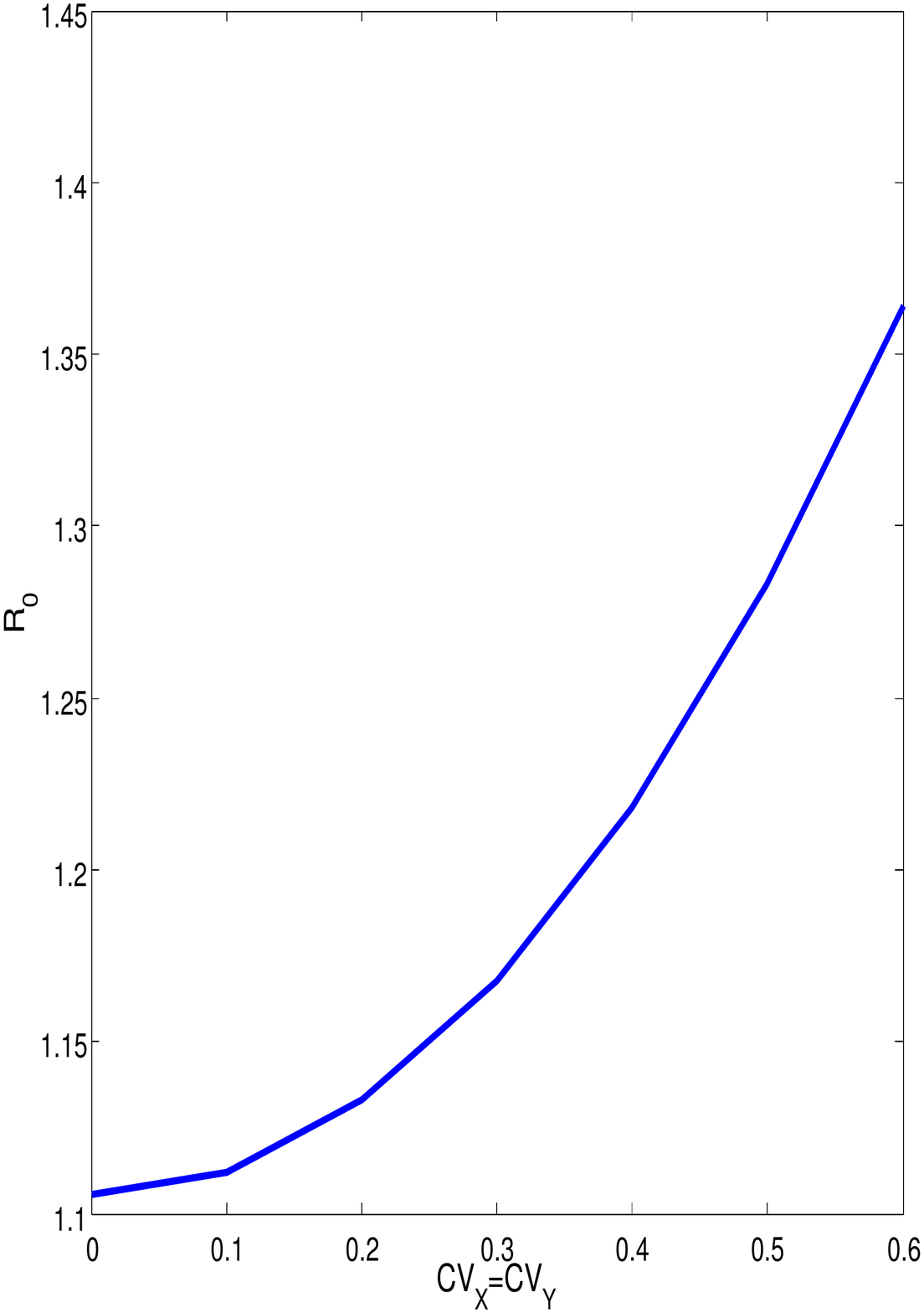}
\caption{$R_0$ as function of $CV_X=CV_Y$ where the weights are random number with a two-point distribution. The degree $D$ and weight $W$ are negatively correlated, and $X$ and $Y$ are  positively correlated (see Example 1 for full specification).}
\label{fig31}
\end{figure}

\noindent \textbf{Example 1}. Our first example is where there is negative correlation between degree and weight (a likely scenario if degree refers to number of sex-partners and weight refers to number of sex-contacts per partner), and where susceptibility and infectivity are positively correlated.
We assume that susceptibility and infectivity are random, each with two possible values, both having mean $0.5$, such that the correlation between susceptibility and infectivity equals $0.8$. Furthermore, the degree distribution follows a Poisson distribution with mean $4$ truncated at $15$. The weight is random with two possible values, $1$ and $10$, with $q(w=1|d)=1-d^{-2}=1-q(w=10|d)$, i.e. larger degree results in higher probability for the smaller weight. In Figure \ref{fig31}, this situation is illustrated. We see that $R_0$ is increasing with the coefficient of variation of infectivity and susceptibility. So, if the degree and weight distributions are negatively correlated, which in a sense homogenizes the community, then introducing heterogeneity in susceptibility and infectivity makes $R_0$ increase.\\

\begin{figure}[!h]
\centering
\includegraphics[width=7cm, height=7cm]{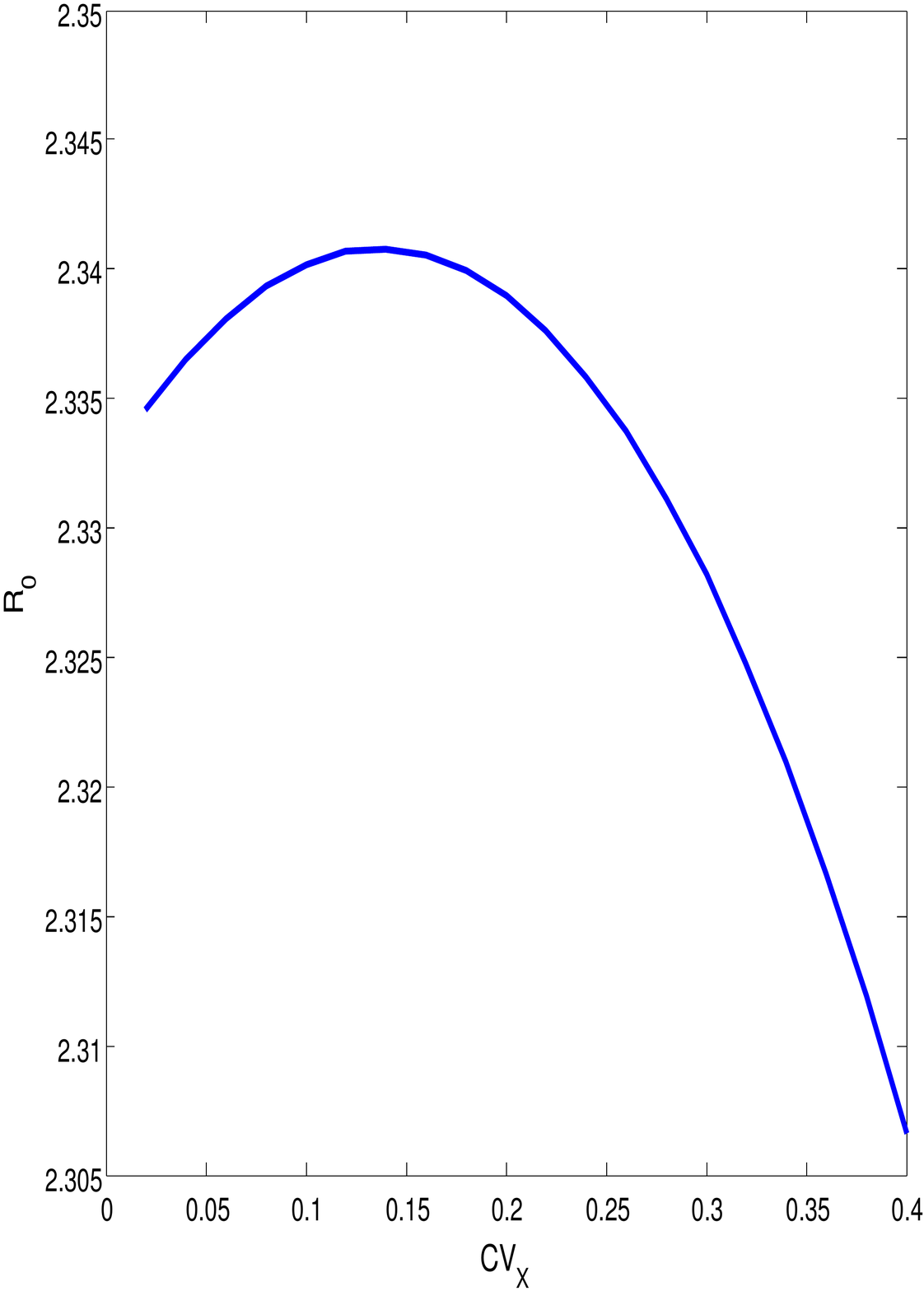}
\caption{$R_0$ as function of $CV_X$, $CV_Y=0.3$, where the weights are random number with two point distribution. The degree $D$ and weight $W$ are uncorrelated, and $X$ and $Y$ are  positively correlated (see Example 2 for full specification).}
\label{fig51}
\end{figure}

\noindent  \textbf{Example 2}. Our second example is where degrees and weights are independent (and keeping positive correlation between susceptibility and infectivity). Just like in Example 1 we let the degree follow a Poisson distribution with mean $4$ truncated at $15$. Also as in Example 1, susceptibility and infectivity each has two possible values, both having mean $0.5$, and correlation $0.8$. However, in order to achieve a non monotone function of $R_0$, we let the coefficient of variation of the infectivity be fixed, $CV_Y=0.3$. The weight is random with two possible values, $1$ and $10$, with $q(1)=q(10)=0.5$, i.e.\ independent of degree. In Figure \ref{fig51} $R_0$ is plotted as a function of $CV_X$ for this situation. We see that $R_0$ is not a monotone function of the coefficient of variation of the susceptibility, but increases initially and at one point start to decrease. This is in contrast to the situation treated in Example 1 where $R_0$ increased monotonically with $CV_X$.\\

\begin{figure}[!h]
\centering
\includegraphics[width=7cm, height=7cm]{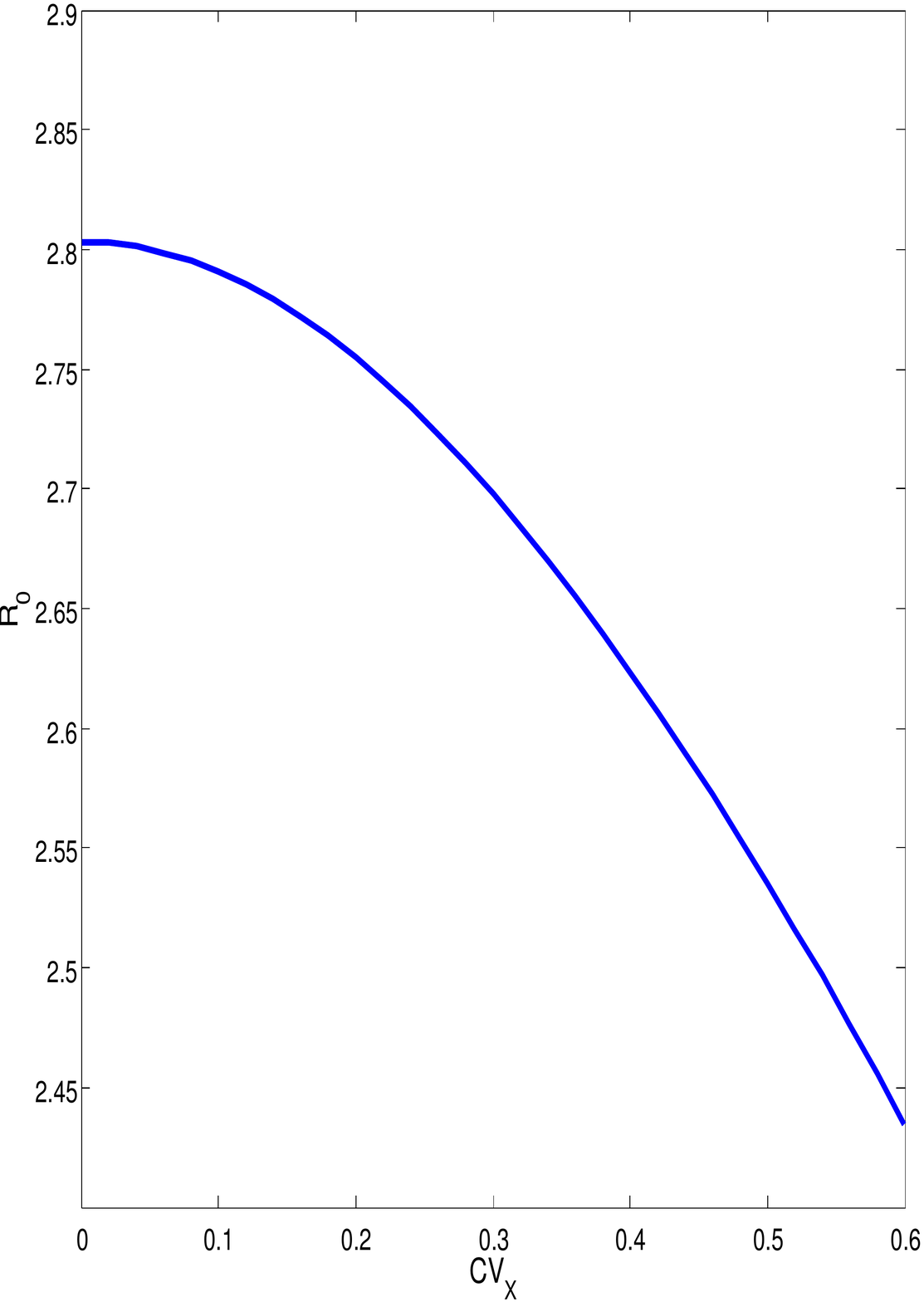}
\caption{$R_0$ as function of $CV_X=CV_Y$ where the weights are random number with two point distribution. The degree $D$ and weight $W$ are positively correlated, also $X$ and $Y$ are positively correlated (see Example 3 for full specification).}
\label{fig21}
\end{figure}

\noindent \textbf{Example 3}. In our last example of this section we look at the case where the degree and weight are positively correlated (keeping positive correlation between susceptibility and infectivity). Susceptibility and infectivity, as well as degree follow the same distribution as in Example 1. The weight is random with two possible values, $1$ and $10$, with $q(w=1|d)=d^{-2}=1-q(w=10|d)$, i.e. now a larger degree results in higher probability of the large weight. This situation is illustrated in Figure \ref{fig21}. We see that $R_0$ now decreases with the coefficient of variation of infectivity and susceptibility. So, if individuals are already heterogeneous in terms of spreading the disease (degree and weight positively correlated), then introducing (independent) heterogeneity in susceptibility and infectivity actually reduces $R_0$.\\

%
%
%

\section{Outbreak probability and size of outbreak}\label{secpi}
The basic reproduction number $R_0$ is not the only informative quantity of an epidemic outbreak. The final proportion infected, $\tau$, and the outbreak probability, $\pi$, are also often of interest. For many epidemic models, including the current model, the initial phase of an outbreak may be approximated by a branching process. Using this approximation, the probability of a major outbreak coincides with the probability of the branching process growing beyond all limits, and this latter probability is well-known how to derive, hence giving a recipe how to compute $\pi$. If, in the epidemic model, the probability for $i$ to infect $j$ is identical to the probability for $j$ to infect $i$, it suffices to have one (undirected) edge present or not between each pair of neighbours, with the interpretation that the edge is present implies that if either gets infected so will the other. The resulting network of undirected edges then specifies who will get infected; those connected to the index case make out the final outbreak. Because in a large community there will be exactly one giant connected component (if $R_0>1$, and no giant if $R_0\le 1$) the probability $\pi$ that randomly selected index case starts a major outbreak is identical to the probability that it belongs to the giant connected component, but the latter is simply the relative size of the giant connected component $\tau$. To conclude we have the well-known property that, if the transmission network can be modeled as an undirected network, then $\pi=\tau$, i.e.\ the branching process approximation for the outbreak probability also gives the relative size $\tau$ of a major outbreak (e.g.\ Britton, 2010).\\

\noindent \textbf{Example 4}. Consider an unweighted network ($W\equiv 1$) with all individuals having the same degree $d$ ($D\equiv d$), and assume that individuals are of two types, with respect to susceptibility and infectivity, at equal proportions.  Assume that both types have susceptibility equal to infectivity, being  $x_1=y_1=\mu_1=\mu-\delta$ and $x_2=y_2=\mu_2=\mu + \delta$, for type 1 and type 2 respectively, with $0 \leq \delta \leq \max(\mu, 1-\mu)$. This implies that $t(w,y_1,x_2)=t(w,x_1,x_2)=t(w,x_2,x_1)=t(w,y_2,x_1)$. Thus the transmission network is undirected in this case. Note that individuals have expected susceptibility and infectivity equal to $\mu$. The offspring matrix $\{m_{i,j}\}$ is defined by $m_{i,j}=(d-1)\mu_i\mu_j/2$ since an infective $i$-individual on average has $(d-1)/2$ susceptible type $j$ neighbours, each of which it infects with probability $\mu_i\mu_j$ ($i=1,2,\ j=1,2$). Simple algebra shows that the coefficient of variation of susceptibility/infectivity is $CV=\delta/\mu$ and $R_0=\mu^2(1+CV^2)(d-1)$, which agrees with Equation (\ref{r0unw2}). Since the probability that the initial infective is of either type with equal probability $0.5$, the probability of a large outbreak is $\pi=(\pi_1+\pi_2)/2$, where $\pi_i$ is the probability of a large outbreak starting with an infective of type $i$. The probability that an $i$-individual infects $k$ type 1-individuals and $l$ type 2-individuals is $$P_i(d,k,l)={d \choose k,l,d-k-l}(\mu_i\mu_1/2)^k(\mu_i\mu_2/2)^l(1-\mu_i\mu_1/2-\mu_i\mu_2/2)^{d-k-l}.$$ From the theory of branching processes, we know that the probabilities $\pi_1$ and $\pi_2$ are given by
\begin{eqnarray}\label{pi22}
1-\pi_1 &=& \sum_{k+l\leq d} (1-\pi'_1)^k(1-\pi'_2)^l  P_1(d,k,l),\label{pi21} \\
1-\pi_2 &=& \sum_{k+l\leq d} (1-\pi'_1)^k(1-\pi'_2)^l P_2(d,k,l),\label{pi22}
\end{eqnarray}
where  $\pi'_1$ and $\pi'_2$ are the solutions to the following system of equations.
\begin{eqnarray}
1-\pi'_1 &=& \sum_{k+l\leq d-1} (1-\pi'_1)^k(1-\pi'_2)^l P_1(d-1,k,l)\label{pi23} \\
1-\pi'_2 &=& \sum_{k+l\leq d-1} (1-\pi'_1)^k(1-\pi'_2)^l P_2(d-1,k,l)\label{pi24}.
\end{eqnarray}
The intuition behind Equations (\ref{pi23}-\ref{pi24}) is the following. To escape a large outbreak starting with one infected individual, none of the offspring may start a large outbreak. Thus, if an infected individual infects $k$ type 1 individuals, and $l$ type 2 individuals, of the $d-1$ susceptible neighbours, no outbreak will occur with probability $(1-\pi'_1)^k(1-\pi'_2)^l$. Taking the expected value w.r.t. the number infected of each type give Equations (\ref{pi23}-\ref{pi24}). Since the index case has $d$ susceptible neighbours (as opposed to $d-1$ in later generations) Equations (\ref{pi21}-\ref{pi22}) follows.\\

\begin{figure}[!h]
\centering
\includegraphics[width=8cm, height=8cm]{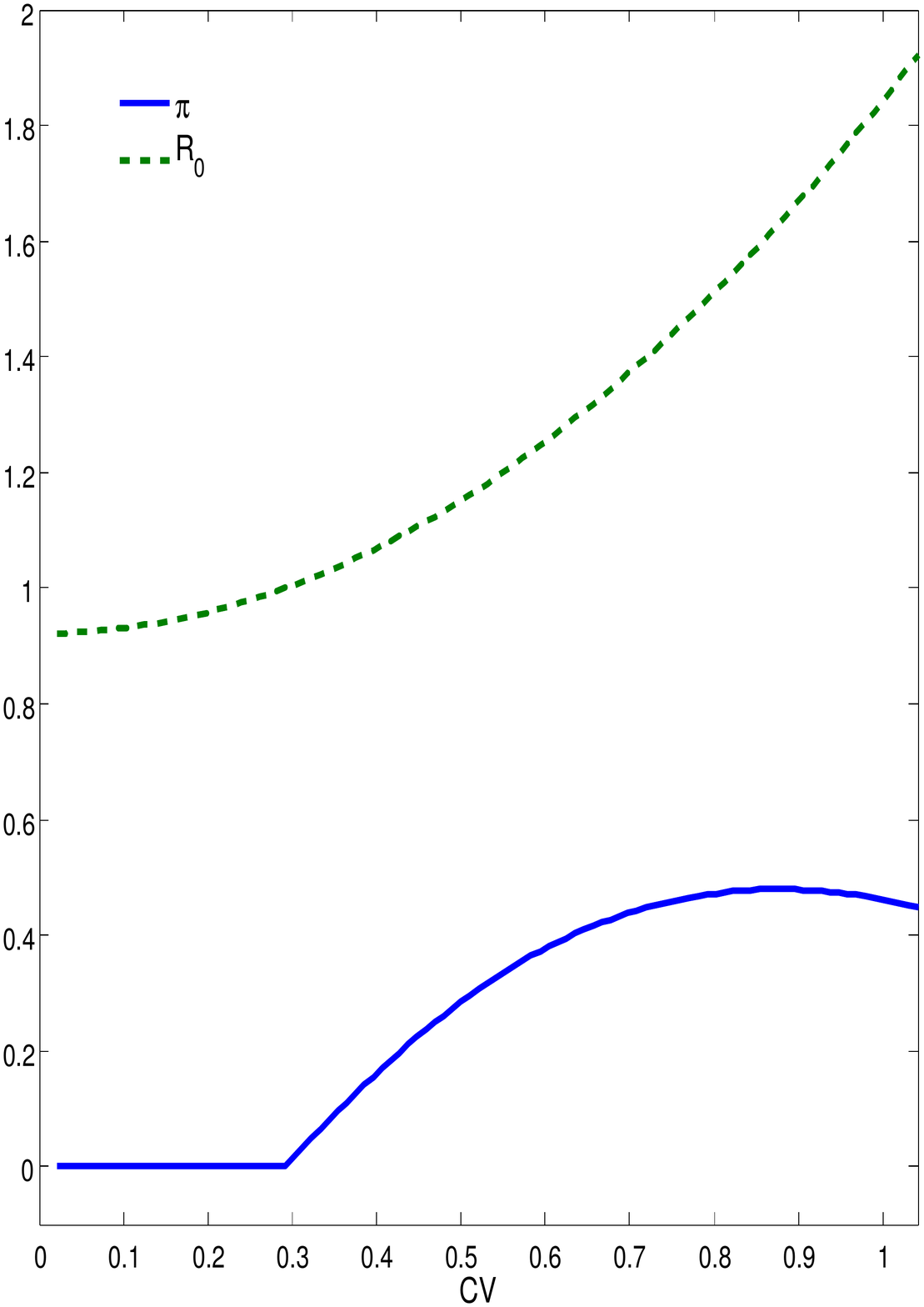}
\caption{The outbreak probability $\pi$, and the basic reproduction number $R_0$, as function of the coefficient of variation of the infectiousness/susceptibility $CV_X=CV_Y$ (see Example 4 for full specification). Since the underlying transmission network is undirected $\pi=\tau$, so the lower curve is hance also a curve of the relative size of a major outbreak if it takes place.}
\label{fig7}
\end{figure}

\noindent This situation is illustrated in Figure \ref{fig7} with $d=5$ and $\mu=0.48$. We see that although $R_0$ is monotonically increasing with the coefficient of variation of susceptibility/infectivity, the outbreak probability $\pi$ is not monotone. Note that $\pi=0$ when $R_0\leq 1$ as to be expected. In the current example we had the transmission probability from one individual to another being the same in both direction. As a consequence, our transmission network may be modelled by an undirected network which, as argued for earlier, implies that the outbreak probability $\pi$ is identical to the relative size $\tau$ of a major outbreak. This hence implies that the final size $\tau$ is not monotone in the coefficient of variation as opposed to $R_0$ which \emph{is} monotonically increasing.
\\


\section{Discussion}\label{secdis}
\noindent In the current paper we have analysed an epidemic model taking place on an edge-weighted network, where the weights on the edges affect the transmission probability. For this model we have focused attention on how individual variation: in degree, weights and in particular susceptibility and infectivity, affect properties of the epidemic. Having no variation in the other factors it was seen that variation in degree \emph{increases} the basic reproduction number $R_0$, and variation in susceptibility and infectivity does the same. For the situation with fixed degree, susceptibility and infectivity we get the opposite effect when introducing variation in weights. The reason for this somewhat surprising opposite effect for weight variation is however explained from the parametrization in the model where weights enter as exponents rather than linearly.\\

\noindent An immediate hypothesis would hence be that introducing randomness in one factor when other factors already have variation, would show a similar pattern, i.e.\ that this would increase $R_0$. This does however not hold in general. For example it was seen in Figure \ref{fig21} that $R_0$ decreased with $CV_X=CV_Y$ when the degree and weights were positively correlated, and in Equation (\ref{R_0-fixed}) when susceptibility and infectivity were negatively correlated. The general qualitative conclusion seems to be that $R_0$ increases when introducing heterogeneity in one aspect unless the remaining factors already induce ''severe'' heterogeneity. In a sense, if individuals are already very heterogeneous in certain aspects, introducing additional (independent) heterogeneity may actually have the effect that individuals become \emph{less} heterogeneous (a phenomenon related to what is known as \emph{regression towards the mean}). \\

\noindent It was also observed that the probability of a major outbreak $\pi$, and the size of an outbreak $\tau$ when the transmission probabilities were reciprocal, had a more complicated dependence on the randomness in degree, weight, susceptibility and infectivity. That is, even when $R_0$ was monotone in $CV_X$, the same is not necessarily true for the outbreak probability $\pi$ and the outbreak size $\tau$ (as shown in Example 4).\\

\noindent The current work may be extended in several ways in order to increase realism. For example, the model does not evolve in real time so the concept of length of infectious period and its randomness is missing. Also, if considering an epidemic over a longer time period, then having a time dynamic network would be of interest to study. Another feature of interest would be to study the effect of vaccination, or other preventive measures, on the reproduction number. Further work might also analyse the current model in more detail and/or construct other examples, to further investigate the qualitative properties of the model.


\section*{Acknowledgements}

Tom Britton is grateful to Riksbankens Jubileumsfond (The Bank of Sweden Tercentenary Foundation) for financial support.


\begin{thebibliography}{10}

\bibitem{A99}
Andersson H.
\newblock {Epidemic models and social networks}.
\newblock {\em Math. Sci.} 24(2):128-147, 1999.

\bibitem{BMS97}
Ball F., Mollison D. and Scalia Tomba G.
\newblock {Epidemics in populations with two levels of mixing}.
\newblock {\em Ann. Appl. Prob.} 7(1):46-89, 1997.

\bibitem{BM90}
Becker N.G. and Marschner, I.C.
\newblock {The effect of heterogeneity on the spread of disease.}
\newblock {\em Lect. Notes Biomath.} 86:90-103, 1990.


\bibitem{B10} Britton, T. 
\newblock {Stochastic epidemic models: a survey}. 
\newblock {\em Math. Biosci.} 225: 24-35, 2010.



\bibitem{BD11}
Britton, T. and Deijfen, M. and Liljeros, F.
\newblock {A Weighted Configuration Model and Inhomogeneous Epidemics}.
\newblock {\em Journal of Statistical Physics}. 145(5):1368-1384, 2011.


\bibitem{DH00}
Diekmann, O. and Heesterbeek, J. A. P.
\newblock {\em Mathematical epidemiology of infectious diseases}.
\newblock {Wiley Series in Mathematical and Computational Biology. John Wiley \& Sons Ltd.}, Chichester, 2000.

\bibitem{J75}
Jagers, P.
\newblock {\em Branching processes with biological applications}.
\newblock {Wiley-Interscience [John Wiley \& Sons]}, London, 1975.

\bibitem{hofstad}
van der Hofstad, R.
\newblock {\em Random Graphs and Complex networks}.
\newblock {Lecture notes}, 2011.
\newblock {http://www.win.tue.nl/~rhofstad/NotesRGCN.pdf}

\bibitem{trapman}
Meester R., Trapman P.
\newblock {Bounding basic characteristics of spatial epidemics with a new percolation model}.
\newblock {\em Adv. Appl. Prob.} 42(2):335-347, 2011.

\bibitem{miller}
Miller, J. C.
\newblock {Bounding the size and probability of epidemics on networks}.
\newblock {\em J. Appl. Probab.} 45(2):498-512, 2008.

\bibitem{newman}
Newman, M. E. J.
\newblock {The structure and function of complex networks}.
\newblock {\em SIAM Rev.} 167-256 (electronic), 2003.


\end{thebibliography}

\section{Appendix}\label{secapp}
\noindent Let $c=\frac{p\mu_W}{1-p}$, then it follows from Equation (\ref{r034}) that
\begin{eqnarray*}
R_0      &=& (d-1) \Big(1-\Big(\frac{(1-p)\frac{r}{\mu_W}}{(1-p)\frac{r}{\mu_W}+p}\Big)^{r}\Big)\\
         &=& (d-1)(1-exp\{-r log(1+c/r)\})\\
(R_0)'_r &=& -exp\{-r log(1+c/r)\}\Big(\frac{-d(rlog(1+c/r)}{dr}\Big)\\
         &=& exp\{-rlog(1+c/r)\}\Big(log(1+c/r)-\frac{c/r}{1+c/r}\Big)
\end{eqnarray*}
for all $r>0$. Define $z=c/r$ and $g(z)=log(1+z)-\frac{z}{1+z}$
where $z>0$. It follows that $g'(z)=1/(1+z)-1/(1+z)^2$. Since $g(0)=g'(0)=0$ and $g'(z)>0$ for $z>0$, it follows that $g(z)>0$ when $z>0$. We conclude that $(R_0)'_r>0$ when $r>0$. Thus $R_0$ is an increasing function of $r$.

\end{document}